\newcommand{\E}[0]{{\sf E}}

\documentclass[letterpaper,11pt]{article}

\usepackage{amsthm}
\usepackage{amsmath}

\newtheorem{thm}{Theorem}[section]

\newtheorem{lemma}[thm]{Lemma}
\newtheorem{prop}[thm]{Proposition}

\newtheorem{conj}[thm]{Conjecture}

\newcommand{\beq}[1]{\begin{equation}\label{#1}}
\newcommand{\enq}[0]{\end{equation}}

\newcommand{\bn}[0]{\bigskip\noindent}
\newcommand{\mn}[0]{\medskip\noindent}
\newcommand{\nin}[0]{\noindent}

\newcommand{\sub}[0]{\subseteq}
\newcommand{\sm}[0]{\setminus}
\renewcommand{\dots}[0]{,\ldots,}

\newcommand{\cee}[0]{{\cal C}}

\newcommand{\eee}[0]{{\cal E}}

\newcommand{\h}[0]{{\cal H}}

\newcommand{\T}[0]{{\cal T}}

\newcommand{\ra}[0]{\rightarrow}
\newcommand{\Ra}[0]{\Rightarrow}

\newcommand{\ZZ}[0]{{\bf Z}}

\newcommand{\0}[0]{\emptyset}

\renewcommand{\qed}[0]{\begin{flushright} \rule{2mm}{3mm} \end{flushright}}

\newcommand{\C}[2]{{{#1}\choose{{#2}}}}
\newcommand{\Cc}[0]{\tbinom}
\newcommand{\ga}[0]{\alpha }
\newcommand{\gb}[0]{\beta }

\newcommand{\gd}[0]{\delta }
\newcommand{\gD}[0]{\Delta }
\newcommand{\gG}[0]{\Gamma }

\newcommand{\gl}[0]{\lambda }

\newcommand{\go}[0]{\omega}
\newcommand{\gO}[0]{\Omega}

\newcommand{\gz}[0]{\zeta}
\newcommand{\eps}[0]{\varepsilon }
\newcommand{\vt}[0]{\vartheta}

\newcommand{\vp}[0]{\varphi}

\newcommand{\sugr}[1]{\mn\textcolor{Red}{[#1]}}

\newcommand{\comments}[1]{}
\usepackage[usenames]{color}
\begin{document}
\renewcommand{\thefootnote}{\fnsymbol{footnote}}

\title{On the triangle space of a random graph\footnotemark}
\author{B. DeMarco, A. Hamm and J. Kahn}
\date{}
\footnotetext{ $^*$Supported by NSF grant DMS0701175.}

\maketitle

\begin{abstract}
Settling a first case of a conjecture of
M. Kahle
on the homology of the clique complex of the random graph $G=G_{n,p}$, we show, roughly speaking, that
(with high probability)
the triangles of $G$ span its cycle space
whenever each of its edges lies in a triangle
(which happens (w.h.p.) when $p$ is at least about $\sqrt{(3/2)\ln n/n}$,
and not below this unless $p$ is very small.)
We give two related proofs of this statement,
together with a
relatively simple proof of a fundamental
``stability" theorem for triangle-free subgraphs of $G_{n,p}$,
originally due to Kohayakawa, \L uczak and R\"odl,
that underlies the first of our proofs.
\end{abstract}

\section{Introduction}\label{Intro}

The primary purpose of this paper is to prove
a first case
(Theorem \ref{TKahle}) of a conjecture
of M. Kahle on the homology of the
clique complex of the (usual)
random graph
$G_{n,p}$.  We will give two (not unrelated) proofs of this.
Underlying the first is
Theorem \ref{8.34}, a (known)
``stability" theorem for
triangle-free subgraphs of $G_{n,p}$, and our second main
contribution is an alternative proof
of this basic result.
We begin with some background.

All graphs will have the vertex set $V=[n]=\{1\dots n\}$,
so we will often fail to distinguish between a graph $G$ and its edge set,
and will tend to regard subgraphs of $G$ as subsets of $E(G)$.
Recall that a {\em cut} of $G$ is $\nabla(W,V\sm W)$, the set of edges
of $G$ joining $W$ and $V\sm W$ for some $W\sub V$.

More or less following \cite{Diestel}, we set, for a given $G$,
$\eee=\eee(G) =(\ZZ/2)^{E(G)}$ (the {\em edge space} of $G$).
We regard elements of $\eee$ as subgraphs of $G$ in the natural way
(namely, identifying a subgraph with its indicator), and write ``$+$"
for symmetric difference.
The {\em cycle space}, $\cee=\cee(G)$ is the subspace of $\eee$
spanned by the cycles, and
$\cee^\perp$ ($:=\{H:\langle F,H\rangle = 0 ~\forall F\in \cee\}$
with the usual inner
product) is precisely the set of cuts (which, note,
includes $\0$).
We are particularly interested in the {\em triangle space},
$\T=\T(G)$, the subspace of $\cee$ spanned by the triangles of $G$.
Recall (see e.g. \cite{Roy})
that the {\em clique complex}, $X(G)$, of a graph $G$ is the
simplicial complex whose faces are the (vertex sets of)
cliques of $G$.

In the rest of this section we write $G$ for $G_{n,p}$,
where, as usual, $p=p(n)$.
A precise possibility, suggested by M. Kahle (\cite{KahlePC}; see also
\cite{Kahle,KahleArxiv}) and proved
by him for $\gG={\bf Q}$ \cite{KahleArxiv}, is
\begin{conj}\label{CKahle}
Let $\gG$ be either $\ZZ$ or a field.
For each positive integer k and $\eps >0$,
if
$$
p>(1+\eps)\left[(1+k/2)(\log n/n)\right]^{1/(k+1)},
$$
then w.h.p. $H_k(X(G),\gG)=0$, where $H_k$ denotes $k$th homology group.
\end{conj}
\nin
(Here $\log $ means $\ln$ and an event holds {\em with high probability}
(w.h.p.) if its probability tends to 1 as $n\ra\infty$.)
We omit topological definitions, since we won't need them in what
follows; see for example \cite{Roy,Kahle}.
For $k=0$---with, of course,
$H$ replaced by the reduced homology $\tilde{H}$---Conjecture \ref{CKahle}
is more or less the classical
result of Erd\H{o}s and R\'enyi \cite{ER1} giving the threshold for
connectivity of $G_{n,p}$.

We will prove Conjecture \ref{CKahle} for
$k=1$ and $\gG=\ZZ/2$, which, being the first unsettled case,
has apparently been the subject of some previous efforts
\cite{Alon,KahlePC}.
Note that here the conclusion ($H_1(X(G),\ZZ/2)=0$) is
just $\T(G) = \cee(G)$, so that the desired statement is
\begin{thm}\label{TKahle}
If $C> \sqrt{3/2}$ is fixed and
$p> C\sqrt{\log n/n}$, then w.h.p. $\T(G)=\cee(G)$.
\end{thm}

\nin
We will actually prove the following more precise version, in
which we set
$
Q= \{\mbox{every edge is in a triangle}\}.
$

\begin{thm}\label{TKahle'}
$$\max_p\Pr(\mbox{$G$ satisfies $Q$ and $\T(G)\neq \cee(G)$}) \ra 0 ~~(n\ra \infty).$$
\end{thm}

\mn
This gives Theorem \ref{TKahle}, since it's
easy to see (see \eqref{EX}) that for $p$ as in that theorem,
w.h.p. every edge of $G$ does lie in a triangle.
Note also (see Proposition \ref{2ndmm})
that for significantly smaller $p$,
$Q$ is unlikely;
so
Theorem \ref{TKahle'} is really about $p$ roughly as in
Theorem \ref{TKahle}.

\medskip
As mentioned above,
we will also give
a new proof of

\begin{thm}\label{8.34}
For each $\eta>0$ there is a $C$ such that
if $p > Cn^{-1/2}$, then w.h.p. each triangle-free
subgraph of $G$ of size at least $|G|/2$
can be made bipartite by deletion of at most $\eta n^2p$ edges.
\end{thm}

\nin
This seminal result---essentially Theorem 8.34 of
\cite{JLR}---seems due to Kohayakawa, \L uczak and
R\"odl \cite{KLR},
though Tomasz \L uczak \cite{LuczakPC} tells us it was
already known (within some small circle) at the time.
(Theorem \ref{8.34} is a slightly restricted version
of the actual result,
corresponding to what's in \cite{JLR};
see Theorem \ref{8.34'} below for the full statement.)

Theorem \ref{8.34} is a ``stability" version of the
following ``density" theorem, which is
essentially
due to Frankl and R\"odl \cite{FR86}.
(More precisely, this is a little stronger than what's stated
in \cite{FR86}, but is easily gotten from
their proof; see also
\cite{HKL} or \cite[Theorem 8.14]{JLR}.)
Write $t(H)$ for the maximum size of a triangle-free subgraph of $H$.

\begin{thm}\label{Dthm}
For each $\gamma>0$ there is a $C$
such that if $p> Cn^{-1/2}$,
then w.h.p.
$t(G) <(1+\gamma)|G|/2$.
\end{thm}
\nin
The relation between Theorems \ref{Dthm} and \ref{8.34}
is like that between Tur\'an's Theorem
\cite{Turan} and the Erd\H{o}s-Simonovits ``stability theorem"
\cite{Simonovits},
which says, roughly, that any $K_r$-free graph with about
$(1-1/(r-1))\C{n}{2} $ edges
is nearly $(r-1)$-partite.
The extension of Theorem \ref{Dthm}
to larger $r$, conjectured in \cite{KLR}, was proved
by Conlon and Gowers
\cite{Conlon-Gowers} and Schacht \cite{Schacht}; the
corresponding extension of
Theorem \ref{8.34},
suggested in \cite{Koh,KLR}, was also proved in \cite{Conlon-Gowers},
then again
in 
\cite{Samotij} (building on \cite{Schacht}), and very recently in
\cite{BMS} and \cite{ST}.
(All of these papers treat more general forbidden subgraphs.)

The original proof of Theorem \ref{8.34}
(see \cite{KLR, JLR}) uses a
sparse version of Szemer\'edi's Regularity Lemma \cite{Szemeredi}
due to Kohayakawa \cite{Koh} and R\"odl (unpublished; see \cite{Koh}),
together with the triangle case of the
``K\L R Conjecture" of \cite{KLR}
(which has
recently been proved in full by Balogh {\em et al.} \cite{BMS}),
while the ingenious recent proofs avoid such tools
(apart from a superficial use of the ``graph removal lemma"
in connection with the present
Lemma \ref{ML3}; see the remark following the statement of
that lemma).

Our (unbiased) feeling is that the argument given
here is the simplest to date,
even compared to specializations of earlier approaches
to the single case covered by Theorem \ref{8.34}
(though Jozsi Balogh \cite{Balogh} tells us
that the specialization of
\cite{BMS} is also
reasonably simple);
it is also of a somewhat different flavor
than earlier work, though there are similarities.
All proofs depend on versions of
Lemma \ref{ML3}.  The argument given here also has in common with
\cite{BMS} and \cite{ST} the use of a small subset of a possible
violator $F$ to significantly restrict the universe from which
the remainder of $F$ must be drawn;
but the mechanism by which we accomplish this is (in a word) more
``dynamic":  it is based on sampling from $G$, while
\cite{BMS} and \cite{ST} depend on an {\em a priori} description
of ({\em all}) triangle-free subgraphs of $K_n$.
The crucial (simple) point supporting our version
is Lemma \ref{ML2}, which seems interesting in itself.

The present proof was obtained independently of \cite{BMS}, \cite{ST}
and was part of the first author's Ph.D. thesis, which was defended around
the time \cite{BMS} and \cite{ST} were posted \cite{DeMarco}.
At this writing we don't know whether the approach can be
extended to prove some of the more general results mentioned above;
generalizing to $K_r$ would just require the corresponding
extension of Lemma \ref{ML2}, which seems true
though we don't yet see a proof.

\medskip
An interest in reproving
Theorem \ref{8.34}---partly motivated by an
application of
that theorem in \cite{DKBPS}---was
actually the starting point for the present work,
as follows.
It's not too hard to show that, roughly speaking,
if $p$ is as in Theorem \ref{8.34}, then w.h.p.
every triangle-free $F\sub G$
with $|F|\geq |G|/2$ has even intersection with most triangles of $G$.
(This again is essentially from \cite{FR86},
following an idea of Goodman \cite{Goodman};
see also \cite[Sec.8.2]{JLR}.)
So in thinking about a new proof of Theorem \ref{8.34},
we wondered whether
some insight might be gained by understanding what happens
when one replaces ``most" by ``all."
This led to the question addressed in Theorem \ref{TKahle},
which we
realized only later was a known problem.

\medskip
The rest of the paper is organized as follows.
Section \ref{LD} consists of various standardish preliminaries,
while Section \ref{MLs}
contains statements of more interesting lemmas,
which are then
proved in Sections \ref{PML1}-\ref{PML3}.
Section \ref{PTK} gives the easy derivation of
Theorem \ref{TKahle'} from Theorem \ref{8.34} and
Lemma \ref{ML1}, and
our proof of Theorem \ref{8.34} is given in
Section \ref{P834}.
Finally, Section \ref{Coda} contains a sketch of a separate proof of
Theorem \ref{TKahle'} that avoids Theorem \ref{8.34}
(this is put off until
the end of the paper to allow reference to
Sections \ref{PTK} and \ref{P834}).

\mn
{\em Usage.}
As noted above, all our graphs will have vertex set
$V=[n]$.
We use
$v\dots z$ for vertices, often without explicitly specifying, e.g.,
``$x\in V$," and $xy$ for the edge more properly written $\{x,y\}$.
We use $|H|=|E(H)|$ (the {\em size} of $H$),
$N_H(x)=\{y:xy\in H\}$ (the {\em neighborhood}
of $x$ in $H$), $d_H(x)=|N_H(x)|$
(the {\em degree} of $x$ in $H$) and $d_H(x,y)= |N_H(x)\cap N_H(y)|$.
For disjoint $S,T\sub V$, $\nabla_H(S,T)$
is the set of edges joining $S,T$
in $H$,
$\nabla_H(S)$ is $\nabla_H(S,V\sm S)$---as noted earlier
such a set of edges,
for which we will often write simply $\Pi$, is a {\em cut} of $H$---and $\nabla_H(v)=\nabla_H(\{v\})$.
As usual, $H[S]$ is the subgraph of $H$ induced by $S$.
We use $T(H)$ for the set of triangles of $H$.

In much of the paper we will take $G=G_{n,p}$ and use this as the
default for $H$, so that (e.g.)
$N(x)=N_G(x)$, $d(x)=d_G(x)$, $\nabla(S,T)=\nabla_G(S,T)$
and, for $B\sub V$,
$N_B(x)=N(x)\cap B$.

Finally, we use $\log$ for $\ln$,
$B(m,p)$ for a random variable with the binomial distribution
${\rm Bin}(m,p)$, and ``$a= (1\pm \vartheta)b$"
for ``$(1-\vt)b\leq a\leq (1+\vt)b$."

\section{Preliminaries}\label{LD}

Here we record some routine probabilistic basics.
We use ``Chernoff's inequality" in the following form (see
\cite[Theorem 2.1]{JLR}), where, for $x\geq -1$,
$\vp(x) = (1+x)\log (1+x)-x$.
\begin{thm}\label{Chern}
For $\xi =B(n,p)$, $\mu=np$ and any $\gl\geq 0$,
$$ 
\Pr(\xi \geq \mu+\gl)\leq \exp [- \tfrac{\gl^2}{2(\mu+\gl/3)}]
$$
and
$$
\Pr(\xi \leq \mu-\gl) \leq \exp[-\mu\vp(-\gl/\mu)]\leq
\exp [- \tfrac{\gl^2}{2\mu}].
$$ 
\end{thm}
\nin
(We will not need the more precise version of the first bound.)

\medskip
We will also (in Section \ref{P834}) need the following
Azuma-Hoeffding type bound.  (A similar statement
can be extracted from, e.g., the discussion in
Section 3 of \cite{aglc} (see (33) and Lemma 3.9(a));
but we include the simple proof.)

\begin{lemma}\label{mart}
Let $X =X(\xi_1\dots \xi_m)$ where the $\xi$'s are i.i.d.,
each with the distribution ${\rm Ber}(p)$, and suppose
$X$ is Lipschitz (that is, changing the value of a single $\xi_i$
changes the value of $X$ by at most 1).
Then for any $t\in [0,1]$, each of
$\Pr(X-\E X < -t)$ and $\Pr(X-\E X > t)$ is at most
$\exp[-t^2/(4mp)]$.
\end{lemma}
\nin
{\em Proof.}
We first observe that if the r.v. $W$ with $\E W=0$ satisfies
$\Pr(W=a)=q = 1-\Pr(W=b)$ for some $a,b$ with $|a-b|\leq 1$,
then for any $\gz\in [0,1]$,
\beq{EegzW}
\E e^{\gz W}\leq e^{-\gz q}[1-q+qe^\gz]\leq e^{\gz^2q},
\enq
where the first inequality follows from the convexity of $e^x$
and the second is an easy Taylor series calculation.

Set $X_i=\E[X|e_1\dots e_i]$, $Z_i= X_i-X_{i-1}$ ($i\in [m]$)
and $Z=\sum Z_i$.
Then
\beq{XEX}
\Pr(X-\E X > t) = \Pr(e^{\gz Z} > e^{\gz t})
<e^{-\gz t}\E e^{\gz Z}.
\enq
while \eqref{EegzW} and induction on $m$
(used in \eqref{mart1} and \eqref{mart2} respectively)
give, again for $\gz\in [0,1]$,
\begin{eqnarray}
\E e^{\gz Z} = \E e^{\gz(Z_1+\cdots + Z_m)}
&=&
\E[\E(e^{\gz(Z_1+\cdots + Z_m)}|\xi_1\dots \xi_{m-1})]
\nonumber\\
&=&
\E[e^{\gz(Z_1+\cdots + Z_{m-1})}\E(e^{\gz Z_m}|\xi_1\dots \xi_{m-1})]
\nonumber\\
&\leq &
\E[e^{\gz(Z_1+\cdots + Z_{m-1})}e^{\gz^2q}]\label{mart1}\\
&\leq &e^{\gz^2 mq}.\label{mart2}
\end{eqnarray}
Finally, inserting this in \eqref{XEX} and taking
$\gz = t/(2mq)$ gives the desired bound.\qed

For the rest of this section we set $G=G_{n,p}$,
and assume
$p$ is at least $n^{-1/2}$.
Of course many of the statements below
hold in more generality,
but there seems no point in worrying about this.
All proofs are quite straightforward, so we give only
one or two representative arguments.

\begin{prop}\label{vdegree}
W.h.p.
\beq{G}
|G| = (1\pm o(1))n^2p/2
\enq
and
\beq{Deg}
d(x)=(1\pm o(1))np ~~\forall x.
\enq
If $p>n^{-1/2}\log^{1/2}n$, then w.h.p.
\beq{CoDeg}
d(x,y)< 4np^2  ~~\forall x,y.
\enq
\end{prop}

\begin{prop}\label{density}
{\rm (a)} For each $\delta$ there is a $K$
such that
w.h.p.
\beq{NST}
|\nabla(S,T)|=(1\pm \delta) |S||T|p
\enq
for all disjoint $S,T\sub V$ of size at least
$Kp^{-1}\log n$.

\mn
{\rm (b)}  For each fixed $\gd>0$, w.h.p.
\beq{NST'}
|\nabla(S)|=(1\pm \delta) |S|(n-|S|)p ~~\forall S\sub V.
\enq
\end{prop}

\medskip
For $X,Y$ (not necessarily disjoint) subsets of $V$, set
$\gz(X,Y)=\gz_G(X,Y)=|\{(x,y)\in X\times Y: xy\in G\}|$.

\begin{prop}\label{gz}
For any $\eps>0$ w.h.p.
\beq{gzYZ}
\gz(Y,Z) = (1\pm \eps)|Y||Z|p
\enq
for all $Y,Z\sub V$ with $|Y||Z|> 8 \eps^{-2}p^{-1}n$.
\end{prop}
\nin
{\em Proof (sketch).}
We may assume $\eps$ is small.
It's easy to see that for a given $Y,Z$,
$\gz(Y,Z)$ can be written as
$B(m_1,p)+B(m_2,p)$ with $m_1+m_2 = |Y||Z|-|Y\cap Z|$.
Failure of \eqref{gzYZ} (for $Y,Z$) then requires that
at least one of these binomials differ from its mean by at least
(essentially) $\eps |Y||Z|p/2$, and the probability of each of these
events is bounded by
$\exp[-\eps^2 |Y||Z|p/[8(1+\eps/3)]$,
which is $o(2^{-n})$ for $Y,Z$ as in the proposition.
\qed

\begin{prop}\label{ST}

\mn
{\rm (a)}
There is a $K$ such that
w.h.p. for all $v$, $S\sub N(v)$ and $T=N(v)\sm S$,
\beq{ST1}
||\nabla (S,T)|- |S||T|p| < Kn^{3/2}p^2
\enq
and
\beq{ST2}
|G[S]| <\left\{\begin{array}{ll}
|S|^2p/2+ Kn^{3/2}p^2&\mbox{in general}\\
o(|S|np^2)&\mbox{if $|S|=o(np)$}
\end{array}\right.
\enq

\mn
{\rm (b)}
There is an $\ga>0$ such that
if $p> 1.2 \sqrt{\log n/n}$ then w.h.p.
\beq{ST3}
|\nabla(S,T)| > \ga|S|np^2
\enq
whenever $v\in V$, $S\sub N(v)$, $T=N(v)\sm S$ and $2\leq |S|\leq |T|$.

\mn
{\rm (c)}
There is a $K$ so that w.h.p.
for all $v$ and $S,T$ disjoint subsets of $N(v)$ with
$|T| > np/3$ and $s>K/p$,
\beq{ST4}|\nabla (S,T)|> 0.9|S||T|p.
\enq
\end{prop}

\mn
{\em Remark.}
The 1.2 in (b) is just a convenient choice between
$1$ and $\sqrt{3/2}$.

\mn
{\em Proof (sketch).}
In each case, by Proposition \ref{vdegree}
(see \eqref{Deg}), it's enough to bound
the probability that the assertion fails at some
$v$ with $d(v) = (1\pm o(1))np$.
We use $s$ and $t$ for $|S|$ and $|T|$.
Having chosen $v$ and $N(v)$ of size
$m=(1\pm o(1))np$, we may bound the number of possibilities
for $(S,T)$ (with given $s,t$) by
$\C{m}{s}$ in (a),(b) and (say)
$\C{m}{s}\C{m}{t}$
in (c).
On the other hand, once we have specified $S $ (and $T$ if we
are in (c)),
we are just bounding a deviation
probability for some binomial random variable, and the required bounds
can (with a little effort) be read off
from Theorem \ref{Chern}.

For example, the most delicate of these assertions is
(b) (which is most delicate for $s=2$).
In general for (b) with $s=o(np)$
(which is more than is needed from (b),
since (a) covers $s $ above about $\sqrt{n}$), we may,
using Theorem \ref{Chern},
bound the probability of a violation with $|S|=s$ by
$$
n\Cc{(1+o(1))np}{s}\exp [- (1-o(1))sn^2p\vp(-1+\ga +o(1))]
~~~~~~~~~~~
$$
\beq{s}
~~~~~~~~~~~~~~~~~~~~ ~~~~~~~~~~~~~~~~~~~~
< n\left\{np\exp[-(1-\gb)np^2]\right\}^s,
\enq
where $\gb=\gb_\ga\ra 0$ as $\ga\ra 0$.
(The initial $n$ is for the choice of $v$,
and we have used
$s=o(np)$ and $d(v) = (1\pm o(1))np$ to say $t= (1-o(1))np$.)
Now
$np\exp[-(1-\gb)np^2]$ is decreasing in $p$, so
is at most $1.2 \sqrt{\log n} ~n^{1/2 - (1-\gb)1.44}$
for $p$ as in (b).
Thus, for slightly small $\ga$, the sum over $s\geq 2$
of the right hand side of \eqref{s} is bounded by
some fixed negative power of $n$.
\qed

\medskip
Finally we should justify the
two comments following the statement of
Theorem \ref{TKahle'}, namely that the property $Q$
(every edge of $G$ is in a triangle) holds w.h.p. if
$p$ is as in Theorem \ref{TKahle} and fails w.h.p. if $p$ is
significantly smaller.  The first of these is trivial:
if $X$ is the number of edges of $G$ not lying in triangles,
then
\beq{EX}
\mu(p):=
\E X = \Cc{n}{2}p(1-p^2)^{n-2},
\enq
which is $o(1)$ for $p > \sqrt{(3/2 +\eps)\log n/n}$
(where, here and in the following proposition, $\eps$ is any
positive constant).
The second assertion is just a
second moment method calculation, whose outcome we record as
\begin{prop}\label{2ndmm}
If $\mu(p) =\go(1)$ then $\Pr(X=0)=o(1)$
(where $X$ and $\mu(p)$ are as in \eqref{EX});
in particular this is true if
$p< \sqrt{(3/2-\eps)\log n/n}$
with $\eps$ a positive constant.
\end{prop}
\nin
\begin{proof}
We have
$X=\sum \E A_{xy}$ with the sum over edges $xy $ of $ K_n$
and $A_{xy}$ the indicator of
$\{\mbox{$xy\in G$ and $xy$ lies in no triangle of $G$}\}$.
We then observe that for $x,y,z,w$ distinct,
$$\E A_{xy}A_{zw} < p^2 (1-p^2)^{2(n-4)}
 ~~~\mbox{and}~~~
\E A_{xy}A_{xz} < p^2 (1-2p^2+p^3)^{n-3},$$
which with \eqref{EX} (and minor calculations which we omit)
gives ${\rm Var(X)}/\E^2X =O(1/\mu(p))$.
\end{proof}

\section{Main lemmas}\label{MLs}

We collect here a few main points underlying
the proofs of Theorems \ref{TKahle'} and \ref{8.34}.
As earlier we write $G$ for $G_{n,p}$.

\medskip
Theorem \ref{TKahle'} says that (for any $p$) it's unlikely
that $Q$ holds but $\T(G)\neq \cee(G)$
(or, equivalently, $\T(G)^\perp\neq \cee(G)^\perp$).
As shown in Section \ref{PTK}, this follows easily from
Theorem \ref{8.34} once we've ruled
out ``small" members of $\T^\perp(G)\sm\cee^\perp(G)$:
\begin{lemma}\label{ML1}
For Q as in Theorem \ref{TKahle'} and fixed $\eta >0$,
\beq{ML1eq}
\max_p\Pr(Q\wedge [\exists F\in \T^\perp(G)\sm\cee^\perp(G),
~|F|< (1-\eta)n^2p/4])<o(1).
\enq
\end{lemma}

\medskip
For a graph $H$ on $[n]$ and $K\sub H$, set
$$B(K,H) = \{e\in K_n\sm H:\mbox{there is no triangle $\{e,f,g\}$
with $f,g\in K$}\}.$$
In the proof of Theorem \ref{8.34}
we will choose $G$ by first choosing
a subgraph $G_0\sim G_{n,\vt p}$ and then placing edges
of $K_n\sm G_0$ in $G$ with probability $(1-\vt)p/(1-\vt p)$
(independently).
Then specification of $F_0=F\cap G_0$, for a triangle-free $F\sub G$,
limits the possibilities for $F\sm G_0$ to subsets of
$B(F_0,G_0)$, and we will want to say this set is small; such an
assertion is supported by
the next lemma (which we will apply
with $G$, $F$ and $p$ replaced by with $G_0$, $F_0$ and $\vt p$).
\begin{lemma}\label{ML2}
For each $\gd>0$ there are $C$ and $\eps >0$ such that
if $p> Cn^{-1/2}$ then
w.h.p.
$|B(F,G)| < (1+\gd) n^2/4 $
for each $F\sub G$ of size at least
$(1-\eps) n^2p/4.$
\end{lemma}

Finally we need the following simple deterministic fact,
in which we write $\tau(F)$ for the number of triangles in $F$.
\begin{lemma}\label{ML3}
If $F\sub K_n$ satisfies $|F|> (1-\gd)n^2/4$
and $|F\sm\Pi|>\eta n^2$ for every cut $\Pi$,
then
$\tau(F)>\tfrac{1}{12}(\eta-3\delta-o(1))n^3$.
\end{lemma}
\nin
{\em Remark.}
As suggested earlier this is not really new, versions of a much
more general statement having
been used in \cite{Conlon-Gowers,Samotij,BMS,ST};
nonetheless we include the simple proof
(see Section \ref{PML1}), both to make
our argument self-contained and to give a reasonable
dependence of $C$ on $\eta$ in Theorem \ref{8.34}.
The results corresponding to Lemma \ref{ML3}
in \cite{Conlon-Gowers,Samotij,BMS,ST} are proved---presumably
just for convenience---using
the ``graph removal lemma" of
\cite{EFR} (so for Lemma \ref{ML3} itself the original
``triangle removal lemma" of Ruzsa and Szemer\'edi
\cite{Ruzsa-Szemeredi}), which for example
gives Lemma \ref{ML3}
with both $\gd$ and
$\tfrac{1}{12}(\eta-3\delta-o(1))$
replaced by some tiny constant
depending on $\eta$.

\section{Proof of Lemma \ref{ML1}}\label{PML1}

We need one easy preliminary observation, which will show up again
in the proof of Theorem \ref{TKahle'}.
\begin{prop}\label{addcut}
Let $G$ be a graph and $F\sub G$, and suppose
$F',F''$ are (respectively) minimum and maximum size members of
$ F + \cee^\perp(G)$.  Then
$$
\forall v
~~ d_{F'}(v)\leq d_{G\sm F'}(v)  ~~\mbox{and} ~~
d_{F''}(v)\geq d_{G\sm F''}(v).
$$
\end{prop}
\nin
(For example if $F'$ violates the first
condition (at $v$), then $F' + \nabla(v)\in F + \cee^\perp(G)$
is smaller than $F'$.)

\medskip
We turn to the proof of Lemma \ref{ML1}, noting that, by Proposition
\ref{2ndmm}, it's enough to bound the probability in
\eqref{ML1eq} when (say) $p > 1.2\sqrt{\log n/n}$,
and for this it's enough to show that the event
in \eqref{ML1eq}---that is,
\beq{ML1eq'}
Q\wedge [\exists F\in \T^\perp(G)\sm\cee^\perp(G),
~|F|< (1-\eta)n^2p/4]
\enq
---cannot occur if $G$ satisfies the conclusions of
Propositions \ref{vdegree}, \ref{density} and \ref{ST}.
Suppose instead that (these conclusions are satisfied and)
\eqref{ML1eq'}
holds, and let $F$ be a smallest member of
$T^\perp(G)\sm\cee^\perp(G)$ and
$J=G\sm F$.  By
Proposition \ref{addcut} we have $d_J(v)\geq d_F(v)$ for all $v$.

For disjoint $S,T\sub V$,
set
$\Psi(S,T)=|\nabla(S,T)|-2|G[S]|$.
Since
$$
\sum_{v} |\nabla(N_F(v),N_J(v))|
=2 |\{T\in T(G): |F\cap T|=2\}|
= 2\sum_{v} |G[N_F(v)]|,
$$
we have
\beq{sumpsi}
\sum_{v} \Psi(N_F(v),N_J(v))=0.
\enq

Let $\eps=\eta/2$
and set
$V_1=\{v:d_F(v)>(1-\eps)np/2\}$,
$V_2= \{v\in V\sm V_1:d_F(v)\geq 2\}$ and
$V_3=V\sm (V_1\cup V_2)$.
Note that $Q$ (with $F\neq \0$, which is all we are now using
from $F\not\in\cee^\perp$)
implies $V_1\cup V_2\neq \0$.
The conclusions of parts (a) and (b) of Proposition \ref{ST}
give, for some fixed positive $\gd$ and $L$,
\begin{eqnarray}\label{sumvpsi}
\sum_{v} \Psi(N_F(v),N_J(v))
 &\geq&
\gd \sum_{v\in V_2}d_F(v) np^2 -L|V_1|n^{3/2}p^2
\nonumber\\
&=&n p^2~[\gd \sum_{v\in V_2}d_F(v) -L|V_1|n^{1/2}].
\end{eqnarray}
(For $v\in V_1$, \eqref{ST1} and \eqref{ST2} give
$$\Psi(N_F(v),N_J(v)) > (d_F(v)d_J(v) -d_F^2(v))p - 3Kn^{3/2}p^2
\geq - 3Kn^{3/2}p^2.$$
A similar discussion gives $\Psi(N_F(v),N_J(v))>\gd d_F(v) np^2$
for $v\in V_2$, where for smaller $d_F(v)$ we use \eqref{ST3}
and the second bound in \eqref{ST2}.)

On the other hand, we will show that
\beq{decent}
\sum_{v\in V_2}d_F(v) =\go(|V_1|n^{1/2}),
\enq
which (with \eqref{sumvpsi})
contradicts \eqref{sumpsi} and completes the proof.

We first observe that \eqref{ST4} implies that
(a.s.) for every $v\in V_1$,
\beq{wNv}
|\{w\in N(v): \min\{|N(w)\cap N_F(v)|,
|N(w)\cap N_J(v)|\} < \tfrac{np^2}{4}\}| < o(np),
\enq
so in particular
\beq{NvV}
|N(v)\cap V_3| = o(np).
\enq
(If $z\in N(v)\cap V_3$, then either
$z\in N_F(v)$, whence $\nabla(z,N_J(v))\sub F$ and
(by the definition of $V_3$) $N(z)\cap N_J(v)=\0$,
or, similarly, $z\in N_J(v)$ and $|N(z)\cap N_F(v)|\leq 1$.)

Now $|F| < (1-\eta)n^2p/4$
implies
$|V_1|< (1-\eps) n$
(since
$(1-\eta)n^2p/4>|F|> (1/2)|V_1|(1-\eps)np/2$
implies $|V_1| < (1-\eta)n/(1-\eps) <(1-\eps)n$).
So by \eqref{NST'}
we have
$$ 
|\nabla(V_1)|> (1-o(1))|V_1|\eps np,
$$ 
which in view of \eqref{NvV} gives
\beq{nabla2}
|\nabla(V_1,V_2)| >(1-o(1))|V_1|\eps np.
\enq
On the other hand, we may assume
$|\nabla_F(V_1,V_2)| =o(|V_1|np)$
(or we have \eqref{decent}),
which gives at least
$(1-o(1))|V_1|\eps np$
pairs $(v,w)$ with
\beq{vw}
\mbox{$v\in V_1$, $w\in V_2$, $vw\in J$ and
$|N_F(w)\cap N_F(v)|> np^2/4$}
\enq
(since by \eqref{wNv} only $o(|V_1|np)$ pairs satisfying
the first three conditions are eliminated by the last).
This gives $\gO(|V_1|np\cdot np^2)$ triples $(v,w,z)$ with $v\in V_1$,
$w\in V_2$, $vw\in J$ and $z\in |N_F(w)\cap N_F(v)|$.
But since each $(w,z)$ belongs to at most $4np^2$ such
triples (see \eqref{CoDeg}), this
says that there are at least $\gO(|V_1|np)$ edges of $F$ meeting $V_2$,
so we have \eqref{decent}.

\section{Proof of Lemma \ref{ML2}}\label{PML2}

\medskip
We prove the lemma with $\eps = .05 \gd$ and $C= 4\eps^{-2}$.
For $F\sub G$, set
$J(F,G) = \{xy\in E(K_n): d_F(x,y)\neq 0\}.$
It is enough to show that for suitable
$C$ and $\eps$, and $p$ as in Lemma \ref{ML2},
w.h.p.
\beq{JFG}
|J(F,G) | > (1-\gd)n^2/4
\enq
for each $F\sub G$ of size at least $(1-\eps)n^2p/4$.
In fact all this needs from the randomization is the property
\beq{gzYZ'}
\gz(Y,Z) = (1\pm \eps)|Y||Z|p
 ~~\forall Y,Z\sub V ~\mbox{with} ~|Y|\geq \eps np ~\mbox{and}~
|Z|\geq \eps n/2,
\enq
which according to Proposition \ref{gz} holds w.h.p.;
thus we assume
\eqref{gzYZ'} holds in $G$ and proceed deterministically.

Given $F\sub G$, set $J=J(F,G)$ and, for
$x\in V$,
$$
\gz(x) = \gz_G(N_F(x),N_J(x)) ~~
(= |\{(y,z):xy\in F, xz\in J,yz\in G\}|).
$$
Then
\beq{N1}
\gz(x) \geq \gz_F((N_F(x),N_J(x))=\sum_{y\in N_F(x)}(d_F(y)-1).
\enq
Heading for a companion upper bound, we
say
$x$ is {\em good} (for $F$) if
$$|\{y\in N_F(x):d_F(y)>\eps np\}|>\eps np$$
(and {\em bad} otherwise),
and let $F^*=\{xy\in F: \mbox{$x,y$ are good}\}$.
We need a few little
observations.
First (we assert)
\beq{F*}
|F\sm F^*|\leq  2\eps n^2p.
\enq
To see this, just notice that
an edge of $F\sm F^*$ either contains a vertex of $F$-degree at most
$\eps np$ or, for some bad $x$, is one of at most $\eps np$ edges of $F$
at $x$ that do not contain a vertex of $F$-degree at most
$\eps np$.

Second, notice that
\beq{dJx}
x ~\mbox{good} ~\Ra ~d_J(x) > \eps n/2.
\enq
For if this fails then there are $Y,Z\sub V$
(namely $Y=N_F(x)$, $Z=N_J(x)$)
with $|Y|>\eps np$, $|Z| \leq \eps n/2$ and
$\gz (Y,Z)\geq |Y|\eps np$,
which implies a violation of \eqref{gzYZ'}
(at $Y$ and some $(\eps n/2)$-superset of $Z$).

Third, again using \eqref{gzYZ'}, we find that if
$x$ is good (or if just $d_F(x)>\eps np$ and
the conclusion of \eqref{dJx} holds) then
$$ 
\gz(x) < (1+\eps) d_F(x)d_J(x)p,
$$ 
which with \eqref{N1} gives (for good $x$)
\begin{eqnarray*}
d_J(x) &> &[(1+\eps)pd_F(x)]^{-1}\sum_{y\in N_F(x)}(d_F(y)-1)\\
&>& \frac{1-\eps}{pd_F(x)}\sum_{y\in N_F(x)}d_F(y),
\end{eqnarray*}
where, since $x$ is good (and $p$ is large),
passing from $(1+\eps)^{-1}$ to $1-\eps$ takes care of the
missing ``$-1$" in the second line.

But then (using \eqref{F*} and our lower bound on $|F|$ in
the last line)
\begin{eqnarray*}
2|J|&\geq &\sum_{x ~\mbox{\small{good}}}d_J(x)
> \frac{1-\eps}{p}\sum_{x ~\mbox{\small{good}}~}
\sum_{y\in N_F(x)}\frac{d_F(y)}{d_F(x)}\\
&\geq &
\frac{1-\eps}{p}\sum_{xy\in F^*}
\left[\frac{d_F(y)}{d_F(x)}+\frac{d_F(x)}{d_F(y)}\right]\\
&\geq & 2(1-\eps)|F^*|/p\\
&>& 2(1-\eps)[(1-\eps)n^2p/4 -2\eps n^2p]/p
>(1-\gd)n^2/2
\end{eqnarray*}
(so we have \eqref{JFG}).

\section{Proof of Lemma \ref{ML3}}\label{PML3}

Suppose $F$ is as in the lemma and
denote by $t_i$ the number of
triangles of $K_n$ containing exactly $i$ edges of $F$,
$i\in \{0,1,2,3\}$ (so $t_3=\tau(F)$).
Writing $X$ for the number of pairs $(e,T)$ with $e\in F$ and $T$
a triangle of $K_n$ containing $e$, we have
\beq{Fn-2}
|F|(n-2)= X=t_1+2t_2+3t_3
\enq
and, according to a nice observation of Goodman \cite{Goodman}
(see \cite[p.209]{JLR} for the easy proof),
\beq{t1t2}
t_1+t_2<n^3/8.
\enq
On the other hand,
\beq{t1t3}
t_1+t_3\geq \eta n^3/3,
\enq
since applying the hypothesized lower bound on
the $|F\sm \Pi|$'s
to the cuts $\Pi=(N_F(v),V\sm N_F(v))$
shows
that each vertex
lies in at least $\eta n^2$
of the triangles counted by $t_1+t_3$.

Now \eqref{Fn-2} and \eqref{t1t2} (together with our assumption on $|F|$)
imply
\begin{align*}
(1-\delta)n^2(n-2)/4<|F|(n-2)&=t_1+2t_2+3t_3\\
&=2(t_1+t_2)-t_1+3t_3<n^3/4-t_1+3t_3,
\end{align*}
whence
\beq{t1t3second}
t_1-3t_3<(\delta+o(1))n^3;
\enq
and combining this with \eqref{t1t3} gives  $t_3>\tfrac{1}{12}(\eta-3\delta-o(1))n^3$.

\section{Proof of Theorem \ref{TKahle'}}\label{PTK}

By Proposition \ref{2ndmm} and Lemma \ref{ML1},
it's enough to show that for
$p > 1.2\sqrt{\log n/n}$ and a fixed
$\eta>0$, it's unlikely that
$\T^\perp(G)$
contains an $F$ for which
\beq{minF}
\min\{|F'|:F'\in F+\cee^\perp(G)\} > (1-\eta)n^2p/4.
\enq
Now if there is such an $F$, then by Proposition \ref{addcut}
there is one of size at least $|G|/2$, and w.h.p. this also
satisfies (say) $|F\sm \nabla(A,B)|> 0.1n^2p$
for each partition $A\cup B$ of $V$;
for, writing $\nabla$ for
$\nabla(A,B)$, we have
\beq{1-eta}
(1-\eta)n^2p/4 < |F+\nabla| = 2|F\sm \nabla| + |\nabla|-|F|
<2|F\sm \nabla| + o(n^2p)  ,
\enq
where we used Proposition \ref{vdegree} (to say $|G|> (1-o(1))n^2p/2$)
and
Proposition \ref{density}(a) (to say $|\nabla|< (1+o(1))n^2p/4$).
But according to Theorem \ref{8.34}, the probability that there is
such an $F$ is $o(1)$ even for $p > Cn^{-1/2}$
(with $C$ as in Theorem \ref{8.34}).

\section{Proof of Theorem \ref{8.34}}\label{P834}

As mentioned in Section \ref{Intro}, we prove
the slightly stronger version from
\cite{KLR}:
\begin{thm}\label{8.34'}
For any $\eta>0$ there are $\eps>0$ and $C$ such that
if $p > Cn^{-1/2}$ then w.h.p. each triangle-free
subgraph of $G$ of size at least $(1-\eps)n^2p/4$
can be made bipartite by deletion of at most $\eta n^2p$ edges.
\end{thm}

\nin
{\em Proof.}
As suggested in Section \ref{MLs},
we
choose $G$ by first choosing
a subgraph $G_0\sim G_{n,\vt p}$ and then placing edges
of $K_n\sm G_0$ in $G_1:=G\sm G_0$ independently, each
with probability $q:=(1-\vt)p/(1-\vt p)$.

Set $\vt = 10^{-5}\eta^2$.
According to
Lemma \ref{ML2}
(and \eqref{G}),
we may choose $\eps<\vt$ and $C$ so that
w.h.p.
\beq{Q1}
|G|\sim n^2p/2,  ~~ ~~ |G_0|\sim n^2\vt p/2
\enq
and
\beq{Q3}
[F_0\sub G_0, ~
|F_0| > (1-2\eps) |G_0|/2 ]~\Ra~
|B(F_0,G_0)| < (1+\vt) n^2/4
\enq
(with $B(\cdot ,\cdot)$ as in Lemma \ref{ML2}).
Let $Q$ be the event that
\eqref{Q1} and \eqref{Q3} occur.

Call $F\sub G$ {\em bad} if it is triangle-free with
$|F|> (1-\eps)n^2p/4$ and $|F\sm \Pi|> \eta n^2p$ for every $\Pi$.
Let $R$ be the event that $G$ contains a bad $F$ and
$S$ ($\sub R$) the event that some $F\sub G$ and $F_0=F\cap G_0$
satisfy
\beq{FF0}
\mbox{$F$ is bad and $|F_0|> (1-2\eps)\vt n^2p/4$.}
\enq
Then
$$\Pr(S|R)\geq \Pr(B((1-\eps)n^2p/4,\vt)> (1-2\eps)\vt n^2p/4) > 1/2$$
(e.g. by Theorem \ref{Chern}, which of course really gives $1-o(1)$
in place of 1/2); so
we will have $\Pr(R) = o(1)$ (which is what we want) if we show
$\Pr(S) = o(1)$, which,
since $Q$ holds w.h.p., is the same as
\beq{QS}
\Pr(Q\wedge S) =o(1).
\enq

Suppose then that $Q$ holds and that $F\sub G$ and
$F_0:=F\cap G_0$
satisfy \eqref{FF0}, and set $F_1=F\sm F_0$
and $B=B(F_0,G_0)$.
By \eqref{Q3} we have
\beq{B}
|B|< (1+\vt)n^2/4.
\enq

Now according to Lemma \ref{ML3}, $B$ must satisfy at least one of
(say)

\mn
(i)  $|B|<(1-0.1\eta)n^2/4$;

\mn
(ii)  there is a cut $\Pi$ for which $|B\sm \Pi|< 0.9\eta n^2$;

\mn
(iii)  $\tau(B) > .04 \eta n^3$.


\medskip
On the other hand, since $F$ is bad (and $F_1\sub G\cap B$), we have:
$$
|G\cap B|\geq |F_1|\geq |F|-|G_0|>(1-3\vt)n^2p/4;
$$
\begin{eqnarray*}
|G\cap (B\sm \Pi)|&=&|(G\cap B)\sm \Pi)|\geq |F_1\sm \Pi|\\
&\geq &|F\sm \Pi|-|G_0|
> (\eta - \vt )n^2p
\end{eqnarray*}
for every cut $\Pi$ of $K_n$;
and
$X:=(G\cap B)\sm F_1$ is a set of edges meeting (i.e. containing
an edge of) each triangle of $G\cap B$, with
$$
|X|=|G\cap B|-|F_1|< |G\cap B| - (1-3\vt)n^2p/4.
$$

Thus if $Q\wedge S$ holds, then there
is an
$F_0\sub G_0$ such that
$B=B(F_0,G_0)$ satisfies \eqref{B}
and one of the following is true:

\mn
(a)   $~|B| < (1-0.1\eta)n^2/4 $ and $ |G\cap B| > (1-3\vt)n^2p/4$;

\mn
(b)
there is a cut $\Pi$ for which
$$\mbox{$|B\sm \Pi| < 0.9\eta n^2$
and $|G\cap (B\sm \Pi)| > (\eta-\vt )n^2p$;}$$

\mn
(c)  $~\tau(B) > .04\eta n^3$ and either
$~|G\cap B| > (1+ .01\eta)n^2p/4$
or there is some $X\sub  G\cap B$
of size at most $.005 \eta n^2p$ meeting all triangles of $G\cap B$.

\medskip
Now---perhaps the main point---if $G_0$ is as in \eqref{Q1}
(much more than
we need here), then the number of possibilities for $F_0$
(once we have chosen $G_0$)
is less than $2^{\vt n^2 p}$.
So for \eqref{QS}
it's enough to show that, for a given $F_0$ (again, with
$B=B(F_0,G_0)$ satisfying \eqref{B}),
each of the events (a)-(c) has probability at most
$o(2^{-\vt n^2p})$.

For (a), (b) and the event $\{|G\cap B| > (1+ .01\eta)n^2p/4\}$ in (c)
this is immediate from
Theorem \ref{Chern}, which bounds the associated
probabilities by expressions $\exp[-f(\eta)n^2p]$,
with the $f(\eta)$'s roughly $.01 \eta^2/8$, $.005\eta$ and
$.0001\eta^2/8$ respectively.
(It may be worth emphasizing that $B$ is determined by $F_0$;
so e.g. in (a) we're interested in the probability
that $G\cap B$ is large given that $B$ is small.
The bound for (b) includes a factor $2^n$ for the number of
possible  $\Pi$'s, which makes no difference since $n^2p=\go(n)$.)

For the second alternative in
(c) it's convenient to speak in terms of the hypergraph $\h$
whose vertices are the edges of $G':=G\cap B$ and whose edges
are the triangles of $G'$.
Let
$e_1\dots e_m$ be the edges of $B$ and let
$Y$ be the minimum size
of a set of edges meeting all triangles of $G'$.
Since $Y$ is a Lipschitz
function of the independent ${\rm Ber}(q)$
indicators ${\bf 1}_{\{e_i\in G'\}}$,
Lemma \ref{mart} gives, for $0\leq t\leq 2mq$,
\beq{Ydev}
\Pr(Y<\E Y  -t) < \exp[-t^2/(4mq)] .
\enq
On the other hand, we will show (assuming
$\tau(B)> .04\eta n^3$ as in (c))
\beq{EY}
\E Y > .01\eta n^2p.
\enq
This will complete the proof,
since \eqref{Ydev} with $t=.005 \eta n^2 p$
(now just using $m< n^2/2$ and $q<p$
and noting that,
for example, $\tau(B) > .04\eta n^3$ implies $t< 2mq$)
then bounds the probability of an $X$ as in (c) by
$\exp[- 10^{-5}\eta^2n^2p] =o(2^{-\vt n^2p})$.

\mn
{\em Proof of} \eqref{EY}.
We actually show the stronger
\beq{Enu*}
\E\nu^*(\h) > .01\eta n^2p,
\enq
where $\nu^*(\h)$
is the fractional matching
number of $\h$ (see e.g. \cite{Lovasz}).
To see this, say a triangle $T$ of $B$ is {\em good}
if it is contained in $G'$ and each of its edges lies in at most
$1.9nq^2$ triangles of $G'$.
Then for any $T\in T(B)$,
\begin{eqnarray}
\Pr(T ~\mbox{is good}) &>& q^3(1-3\Pr(B(n,q^2)>1.9nq^2))
\nonumber\\
&>& q^3(1-3\exp[-nq^2/4]).
\label{Prgood}
\end{eqnarray}
Define a (random) weighting $w$ of the triangles of $G'$ by
$$
w(T) =\left\{\begin{array}{cl}
(1.9nq^2)^{-1}&\mbox{if $T$ is good,}\\
0&\mbox{otherwise.}
\end{array}\right.
$$
Then $w$ is a fractional matching of $\h$, and we have
(using \eqref{Prgood})
$$
\E \nu^*(\h) \geq \tau(B) (1-3\exp[-nq^2/4])q^3(1.9nq^2)^{-1}
> .01\eta n^2p.
$$
\qed

\section{Coda}
\label{Coda}

Here we sketch an alternate
proof of
Theorem \ref{TKahle}.
The argument is similar to that in Sections
\ref{PTK} and \ref{P834}, but
seems worth including, as it is a little easier and shows that
Theorem \ref{8.34} and Lemma \ref{ML2} were not really
needed.

\medskip
As in Section \ref{PTK},
we just need to show that for
$p > 1.2\sqrt{\log n/n}$ and a fixed (small)
$\eta>0$, it's unlikely that
$\T^\perp(G)$
contains an $F$ for which
\eqref{minF} holds.
We again fix some small $\vt $
(e.g. $\vt=0.1\eta^2$) and
choose $G$ by first choosing
$G_0\sim G_{n,\vt p}$
and then adding edges
of $K_n\sm G_0$
with probability
$(1-\vt)p/(1-\vt p)$.
Of course we again have \eqref{Q1} w.h.p., and a discussion like
that for \eqref{1-eta} shows that w.h.p.
any $F$ with \eqref{minF} satisfies
$|F\sm \Pi|> 0.1n^2p$
for every cut $\Pi$.

Given $G_0$ and $F_0\sub G_0$, set
$A(G_0)=\{xy\in K_n\sm G_0: N_{G_0}(x,y)\neq\0\}$,
$J=J(G_0)  = K_n\sm (G_0\cup A(G_0))$,
and
$$B=B(F_0,G_0)=\{xy\in A(G_0):z\in  N_{G_0}(x,y)\Ra
|\{xz,yz\}\cap F_0|=1\}.$$
Then any $F\in \T^\perp(G)$ with $F\cap G_0=F_0$ satisfies
$
F\sm (F_0\cup J) = G\cap B.
$

Note also that w.h.p.
\beq{GJ}
|G\cap J|< o(n^2p)
\enq
(e.g. by Theorem \ref{Chern}, using
$\E |J| < n^2(1-(\vt p)^2)^{n-2}$ and Markov's Inequality
to say that w.h.p. $|J|=o(n^2)$).

But if \eqref{Q1} and \eqref{GJ} hold and
$F\in \T^\perp(G)$
satisfies \eqref{minF}, then there is an $F_0\sub G_0$ such that,
with notation as above, we have one of:

\mn
(a)  $|B|< (1-2\eta)n^2/4$ and $|G\cap B| > (1-\eta - 2\vt -o(1))n^2p/4$;

\mn
(b)  there is a cut $\Pi$ with
$$\mbox{$|B\sm \Pi| < 0.05 n^2~$
and $~|G\cap (B\sm \Pi)| > (0.1-\vt/2-o(1))n^2p$;}$$

\mn
(c)  $~\tau(B) > .004 n^3$ and $G\cap B$ is triangle-free.

\mn
Here we used \eqref{minF}
(to say $|F|>(1-\eta)n^2p/4$),
\eqref{Q1} and \eqref{GJ} for the second bound in (a);
$|G\cap (B\sm \Pi)|= |F\sm (\Pi\cup F_0\cup J)|$ together
with $|F\sm \Pi|> 0.1n^2 p$, \eqref{Q1} and \eqref{GJ}
for the second bound in (b);
and Lemma \ref{ML3} to say that failure of the
conditions on $B$ in (a) and (b) implies the one in (c).

But, as in Section \ref{P834},
\eqref{Q1} bounds the number of possibilities for $F_0$
(given $G_0$)
by $2^{\vt n^2 p}$, whereas, we assert, each of the events in (a)-(c)
has probability $o(2^{-\vt n^2 p})$.
For (a) and (b) this is again given by Theorem \ref{Chern}
(with (a) dictating the above choice of $\vt$).
For (c) we may, for example, use an inequality of Janson
(\cite{Janson}; see also \cite[Theorem 2.14]{JLR}), as follows.
Write $S$ for the set of (edge sets of) triangles of $B$;
for $A\in S$,
let $I_A$ be the indicator of $\{A\sub G\}$;
and set $m=\tau(B)> .004n^3$.
Then $\mu:= \sum\E I_A=mp^3$ and
$\bar{\gD}:= \sum \sum_{A\cap B\neq \0}\E I_AI_B < 3mnp^5+\mu$,
and Janson's inequality bounds
the probability that $G\cap B$ is triangle-free by (say)
$$\exp [-\mu^2/(2\bar{\gD})] <\exp[-.0006 n^2p].$$\qed

\bn
Department of Mathematics\\
Rutgers University\\
Piscataway NJ 08854\\
rdemarco@math.rutgers.edu\\
hammac3@math.rutgers.edu\\
jkahn@math.rutgers.edu

\end{document}